\documentclass[a4paper,11pt,leqno]{amsart}

\usepackage{amsmath,amssymb,array,euscript,amsfonts}
\usepackage{graphicx,epsfig,color}
\usepackage[usenames,dvipsnames]{xcolor}
\usepackage{rotating}
\usepackage{bbm}
\usepackage{color}
\usepackage{hyperref}

\def\R{\hbox{\font\dubl=msbm10 scaled 1100 {\dubl R}}}

\def\N{\hbox{\font\dubl=msbm10 scaled 1100 {\dubl N}}}

\def\d{\,{\rm{d}}}

\newtheorem{Theorem}{Theorem}
\newtheorem{Lemma}{Lemma}
\title[The $a$-points of the Riemann zeta-function]{The $a$-points of the Riemann zeta-function and the functional equation}

\thanks{\bf Dedicated to the memory of Professor Eduard Wirsing}

\author[Athanasios Sourmelidis, J\"orn Steuding, and Ade Irma Suriajaya]{Athanasios Sourmelidis, J\"orn Steuding, and Ade Irma Suriajaya}

\date{April 2022}

\begin{document}

\begin{abstract}
We prove an equivalent of the Riemann hypothesis in terms of the functional equation (in its asymmetrical form) and the $a$-points of the zeta-function, i.e., the roots of the equation $\zeta(s)=a$, where $a$ is an arbitrary fixed complex number. \end{abstract}

\maketitle

{\small \noindent {\sc Keywords:} Riemann zeta-function, Riemann hypothesis, $a$-points, functional equation}

\section{Motivation and Statement of the Main Results}

Let $s=\sigma+it$ be a complex variable. The Riemann zeta-function $\zeta$ is for $\sigma>1$ defined by
$$
\zeta(s)=\prod_p(1-p^{-s})^{-1},
$$
where the product is taken over all prime numbers $p$, and by analytic continuation elsewhere except for a simple pole at $s=1$.  The Euler product representation above indicates the relevance of $\zeta$ for the distribution of prime numbers. The yet unproven Riemann hypothesis claims that all nontrivial (non-real) zeros $\rho=\beta+i\gamma$ of the Riemann zeta-function $\zeta$ lie on the critical line $1/2+i\R$. This distribution of zeros would lead to the least possible error term in the prime number theorem. 

We define the Chebyshev function by
$$
\psi(x)=\sum_{n\leq x}\Lambda(n),
$$
where $\Lambda$ denotes the von Mangoldt-function that counts prime powers $n=p^k$ with logarithmic weight $\log p$. The prime number theorem is the asymptotic formula $\psi(x)\sim x$, and the Riemann hypothesis is equivalent to 
\begin{equation}\label{rh}
\psi(x)=x+O_\epsilon(x^{1/2+\epsilon}),
\end{equation}
where here and in the sequel $\epsilon>0$ is arbitrary; this equivalence was first proved by Helge von Koch \cite{vkoch} and relies on the explicit formula from Bernhard Riemann's path-breaking memoir \cite{rie}.

Another more analytic aspect is the distribution of values of $\zeta$. Given a complex number $a$, the roots of the equation $\zeta(s)=a$ are called $a$-points; we denote these roots in the right half-plane by $\rho_a=\beta_a+i\gamma_a$ and their count is very similar to the number of nontrivial zeros (i.e., the case $a=0$). It was Edmund Landau who suggested in his invited lecture \cite{ede} at the International Congress of Mathematicians in Cambridge 1912 to study the distribution of $a$-points with the words: ``{\it The points at which an analytic function is equal to $0$ are very important; but equally interesting are the points at which it assumes a certain value $a$.}''\footnote{The authors' translation of the German original text: ``{\it Es sind be einer analytischen Funktion die Punkte, an denen sie $0$ ist, zwar sehr wichtig; ebenso interessant sind aber die Punkte, an denen sie einen bestimmten Wert $a$ annimmt.}''}. The so-called Nevanlinna theory, or value distribution theory, developed by Rolf Nevanlinna a little later, takes up this idea. 

In this note we investigate the $a$-points of $\zeta$ in the context of the functional equation, that is
\begin{equation}\label{feq}
\zeta(s)=\Delta(s)\zeta(1-s),
\end{equation}
where 
\begin{equation}\label{elt}
\Delta(s):=2(2\pi)^{-s}\sin{\textstyle{{\pi s\over 2}}}\,\Gamma(1-s).
\end{equation}
Note that we could as well consider other functions satisfying a similar functional equation, e.g. Dirichlet $L$-functions. However, for the sake of simplicity, we restrict to the case of $\zeta$. 

Our first result deals with an asymptotic formula for $\Delta$ at the $a$-points. 

\begin{Theorem}\label{2}
Let $a$ be an arbitrary fixed complex number. Then, as $T\to\infty$,
$$
\sum_{\beta_a\geq 0,\, 0<\gamma_a<T}\Delta(\rho_a)=a{T\over 2\pi}\log{T\over 2\pi e}-\psi\left({T\over 2\pi}\right) + O_\epsilon\left(T^{1/2+\epsilon}\right),
$$
Moreover, the Riemann hypothesis is true if, and only if, 
$$
\sum_{\beta_a\geq 0,\, 0<\gamma_a<T}\Delta(\rho_a)=a{T\over 2\pi}\log{T\over 2\pi e}-{T\over 2\pi}+O_\epsilon\left(T^{1/2+\epsilon}\right).
$$
\end{Theorem}

\noindent 
All implicit constants in the error terms here and below depend on $a$. Here and in the sequel we consider only $a$-points in the right half-plane; in what follows, however, we do not write explicitly $\beta_a\geq 0$ anymore. Note that all but finitely many $a$-points in the left half-plane lie arbitrarily close to the trivial (real) $\zeta$-zeros (as an application of Rouch\'e's theorem implies).

Using the prime number theorem in combination with the best known zero-free region for $\zeta$ (found by Nikolay Korobov and (independently) Ivan Vinogradov; cf. \cite[Theorem 12.2]{ivic}), leads to the unconditional asymptotic formula
$$
\sum_{\beta_a\geq 0,\, 0<\gamma_a<T}\Delta(\rho_a)=a{T\over 2\pi}\log{T\over 2\pi e}-{T\over 2\pi}+O\big(T\exp\big(-c(\log T)^{3/5}(\log\log T)^{-1/5}\big)\big),
$$ 
where $c>0$ is an absolute constant, independent of $a$. Note that Eduard Wirsing \cite{wirsing} gave an elementary proof of the prime number theorem with remainder term of order $x/(\log x)^c$, where $c>0$ is arbitrary; the same result was achieved around the same time by Enrico Bombieri \cite{bombi}. Although this is remarkable in terms of the used machinery, the result falls short of what analytic methods can achieve.  

The number $N(T)$ of nontrivial $\zeta$-zeros $\rho=\beta+i\gamma$ satisfying $0<\gamma<T$ (counting multiplicities) is ruled by the von Mangoldt-formula, that is
\begin{equation}\label{rvm}
N(T)={T\over 2\pi}\log {T\over 2\pi e}+O(\log T).
\end{equation}
For the number $N_a(T)$ of $a$-points $\rho_a=\beta_a+i\gamma_a$ with $0<\gamma_a<T$ in the right half-plane (counting multiplicities) we have as well 
\begin{eqnarray}\label{rvma}
N_a(T)= {T\over 2\pi}\log{T\over 2\pi ec_a}+O(\log T)
\sim {T\over 2\pi}\log T \sim N(T),
\end{eqnarray}
where $c_a=1$ for $a\neq 1$ and $c_1=2$, which was first shown by Landau \cite{bohr}.

Hence, by Theorem \ref{2}, for any complex number $a$, the mean-value of $\Delta$ at the $a$-points of $\zeta$ in the right half-plane exists and 
$$
\lim_{T\to\infty}{1\over N_a(T)}\sum_{0<\gamma_a<T}\Delta(\rho_a)=a.
$$
This limit is minimal in absolute value for $a=0$, the case of $\zeta$-zeros. 

This special case had been considered earlier (as the authors learned lately). It will be shown in the proof of Theorem \ref{2} that, as $T\to\infty$, we have unconditionally 
$$
\psi(x)=-\sum_{0<\gamma<T}\Delta(\rho)+O_\epsilon\big(T^{1/2+\epsilon}\big),
$$
where $2\pi x=T$, from which we immediately deduce via (\ref{rh}) that the Riemann hypothesis is true if, and only if,
$$
\sum_{0<\gamma<T}\Delta(\rho)=-{T\over 2\pi}+O_\epsilon\big(T^{1/2+\epsilon}\big).
$$
This result is essentially included in a paper of Johannes Schoi\ss{}engeier \cite{schoiss} (by a different method) and, later, Juan Arias de Reyna \cite{arias} (by a similar reasoning).

Concerning the reflected arguments $1-\rho_a$ in place of the $a$-points, we have a rather different picture:

\begin{Theorem}\label{1}
Let $a$ be a complex number $\neq 0$. Then, as $T\to\infty$,
$$
\sum_{0<\gamma_a<T}\Delta(1-\rho_a)={1\over a}\sum_{0<\gamma_a<T}\zeta(1-\rho_a)=\Omega_\epsilon\big(T^{b_a-1/2-\epsilon}\big),
$$
where $b_a$ is the supremum of the real parts $\beta_a$. 
\end{Theorem}

\noindent The first equality follows directly immediately from the functional equation. For the second equality recall that we write $f(t)=\Omega(g(t))$ if $f(t)=o(g(t))$ does not hold as $t\to\infty$.

The case $a=0$ is, of course, different since the nontrivial $\zeta$-zeros are symmetrically distributed with respect to the critical line and the real axis. The case $a=1$ is a little different for technical reasons since then the logarithmic derivative of $\zeta(s)-a$ (which is an important tool in our reasoning below) is not an ordinary Dirichlet series. For other values of $a$, however, one has $b_a>1$, and $b_{|a|}\to\infty$ as $a$ tends to $1$ (as follows from (\ref{aaab}) in the subsequent section). Hence, it follows that the mean-value of $\Delta(1-\rho_a)$ (resp. of $\Delta(\rho_a)^{-1}$) does not exist for uncountably many $a$ (very likely for all $a\neq 0$), while the mean-value of $\Delta(\rho_a)$ always exists (by Theorem \ref{2}). 
\smallskip

In the following section we collect some results that will be needed in the later proofs, namely the proof of Theorem \ref{2} for 
the case $a=0$ in Section 3, resp. the  general case $a\neq 1$ in Section 4, and the proof of Theorem \ref{1} in Section 5.

\section{Preliminaries about $a$-points and Gonek's lemma}

In 1911, Harald Bohr \cite{boh} was the first to prove that $\zeta(s)$ takes any complex value different from zero in the right half-plane $\sigma>1$. Another proof of this result was given by Landau \cite{1933} whose method was based on consideration of the reciprocals of Dirichlet series $f(s)$. Landau proved that if $f(s)$  is an ordinary Dirichlet series with non-zero constant term, then $1/f(s)$ can also be represented as an ordinary Dirichlet series in the zero-free right half-plane of $f(s)$, which always exists (see for example \cite{tit} \S 9.6). In fact, it can be shown that $\sup\left\{\sigma:f(s)=0\right\}$ is the abscissa of convergence of $1/f(s)$.

Therefore, if $a\neq0,1$, then $(\zeta(s)-a)^{-1}$ can be represented as an ordinary Dirichlet series whose abscissa of convergence is $b_a=\sup\beta_a$. Based on the discussion above, we deduce that $1<b_a<\infty$. Moreover, if $\overline{b}_a$ denotes the abscissa of absolute convergence of  $(\zeta(s)-a)^{-1}$, then the product 
$$
\frac{\zeta'(s)}{\zeta(s)-a}=\sum_{n\geq1}\frac{\Lambda_a(n)}{n^s}
$$
is also an ordinary Dirichlet series which is absolutely convergent for $\sigma>\overline{b}_a$. The implicitly defined coefficients $\Lambda_a$ generalize the von Mangoldt-function $\Lambda=-\Lambda_0$. These results and more can be found in Landau's article, as well as for example in the work of Bombieri \& Amit Ghosh \cite[Section 3]{bomgh} and Siegfred Baluyot \& Steven Gonek \cite{balgon}.

In view of the inequalities $0< \vert\zeta(s)\vert \leq \zeta(\sigma)$, valid for $\sigma>1$, and $\zeta(s)=1+2^{-s}(1+o(1))$ as $\sigma\to+\infty$, it follows that the $a$-points cannot lie too far to the right for $a\neq 1$. More precisely, one has
\begin{equation}\label{aaab}
b_{\vert a\vert}\sim {1\over \log 2}\log{1\over \vert a\vert-1},
\end{equation}
as $\vert a\vert\to\infty$. In the special case $a=1$, when $a$ is equal to the constant term of the Dirichlet series for $\zeta$, one considers analogously the Dirichlet series $2^s\sum_{n\geq 2}n^{-s}$. 

The following result about the value-distribution of the logarithm of the zeta-function belongs to Bohr \& B{\o}rge Jessen \cite{boje}. 

\begin{Lemma}\label{boj} 
If $\log\zeta(s)$ comes arbitrarily near to a given number $c$ on a vertical line $\sigma_0+i\R$, where $\sigma_0>1$, then in every strip $\sigma_0-\delta<\sigma<\sigma_0+\delta$  the value $c$ is taken more than $\kappa(c,\sigma_0,\delta)T$ times, for large $T$, in $0<t<T$, where $\kappa(c,\sigma_0,\delta)$ is a positive constant.

Moreover, there exist positive constants $\kappa_1(c)<\kappa_2(c)$, depending on $c$, such that the number $M_c(T)$ of zeros of $\log\zeta(s)-c$ in $\sigma>1$ satisfies 
$$
\kappa_1(c)T<M_c(T)<\kappa_2(c)T.
$$
\end{Lemma}

\noindent For a proof see \cite{boje} and \cite{titch}, \S 11.8, respectively. It is obvious how this applies to the $a$-points of $\zeta$ via taking $a=\exp(c)$.

In addition, assuming the Riemann hypothesis, already Landau \cite{bohr} showed that the $a$-points are clustered around the critical line, and Norman Levinson \cite{levi} improved upon this by proving unconditionally that most of the $a$-points have real part arbitrarily close to $1/2$ (in a quantitative way).  

We need in addition the partial fraction decomposition of the logarithmic derivative of $\zeta(s)-a$, namely
\begin{equation}\label{partfraca}
{\zeta'(s)\over \zeta(s)-a}=\sum_{\vert t-\gamma_a\vert\leq 1}{1\over s-\rho_a}+O(\log (2+\vert t\vert)),
\end{equation}
valid for $|t|\geq 1$ and $-1\leq\sigma\leq\sigma_0$ with any fixed positive real number $\sigma_0$. This is a straightforward generalization of the classical formula for $a=0$; a proof of (\ref{partfraca}) can be found in \cite{gs}, Lemma 8. In view of (\ref{rvm}) and (\ref{rvma}) this implies the bound
\begin{equation}\label{partfrac}
{\zeta'(s)\over \zeta(s)-a}\ll \log (2+\vert t\vert)^2,
\end{equation}
valid for the same range as (\ref{partfraca}).

We conclude with two standard results. First of all, a bound for $\Delta$ (defined by (\ref{elt})) that follows immediately from Stirling's formula, i.e.,
$$
\Delta(\sigma+it)= \left({t\over 2\pi}\right)^{1/2-\sigma}\exp\left(-it\log{t\over 2\pi e}+{\pi i\over 4}\right)\left(1+O(t^{-1})\right)
$$
for $t\geq 1$, and
\begin{equation}\label{delt}
\Delta(\sigma+it) \asymp  \vert t\vert^{1/2-\sigma}
\end{equation}
as $\vert t\vert\to\infty$; see for example \cite{titch}, (4.12.3).

Finally, we need a modified version of a lemma due to Gonek \cite{gonek}\footnote{The additional log factor in the error term does not exist in the original article. However,  having a closer look in the proof shows that it cannot be omitted.}:

\begin{Lemma}\label{gonek}
Let $(b_n)_{n\geq1}$ be a sequence of complex numbers such that $b_n\ll_\epsilon n^{d+\epsilon}$ for some $d\geq0$. If $c>d$ and $m\geq0$ is an integer, then as $T\to\infty$
\begin{eqnarray*}
\lefteqn{\frac{1}{2\pi i}\int_{c+i}^{c+iT}\Delta^{(m)}(1-s)\sum_{n\geq 1}b_n n^{-s}\d s}\\
&=&\sum_{n\leq T/(2\pi)}b_n(\log n)^{m}+O\left(T^{c-1/2}(\log T)^{m+1}\right).
\end{eqnarray*}
\end{Lemma}

\noindent This result is also our starting point for the proof of Theorem \ref{2}.

\section{Proof of the special case $a=0$}

We apply Lemma \ref{gonek} to the logarithmic derivative 
$$
-{\zeta'\over \zeta}(s)=\sum_{n\geq 2}\Lambda(n)n^{-s},
$$
where the Dirichlet series on the right converges absolutely in the half-plane $\sigma>1$ and uniformly in every compact subset. 
This yields
\begin{eqnarray*}
-\psi(x)=-\sum_{n\leq x}\Lambda(n)&=&{1\over 2\pi i}\int_{1+\delta+i}^{1+\delta+iT}\Delta(1-s){\zeta'\over \zeta}(s)\d s+O_\delta(T^{1/2+2\delta}),
\end{eqnarray*}
where $\delta>0$ is fixed and $T=2\pi x$. 

We assume that $T$ is at a distance from any $\zeta$-zero; more precisely, $\vert T-\gamma\vert\gg (\log T)^{-1}$ by (\ref{rvm}); this restriction can later be removed at the expense of an error term of size $O(T^{1/2+\delta}\log^2 T)$ by \eqref{partfrac} and (\ref{delt}). 

Shifting the path of integration to the left, it follows from Cauchy's theorem that
\begin{eqnarray*}
\int_{1+\delta+i}^{1+\delta+iT}\Delta(1-s){\zeta'\over \zeta}(s)\d s=\left\{\int_{1+\delta+i}^{-\delta+i}+\int_{-\delta+i}^{-\delta+iT}+\int_{-\delta+iT}^{1+\delta+iT} \right\}\ldots +2\pi i\,\Sigma,
\end{eqnarray*}
where $\Sigma$ is the sum of residues at the nontrivial $\zeta$-zeros in the interior of the rectangular contour; here we suppose that $\delta\in(0,1)$ in order to avoid a contribution from trivial zeros. Notice that $\Delta(s)$ is regular except for simple poles at $s=1+2n$ for $n\in\N$. 

The lower horizontal integral is bounded by a constant. For estimating the other integrals we apply in a straightforward manner (\ref{delt}) as well as (\ref{partfrac}) in combination with (\ref{rvm}), and get 
\begin{eqnarray*}
\int_{-\delta+i}^{-\delta+iT}\Delta(1-s){\zeta'\over \zeta}(s)\d s\ll \int_1^T t^{-1/2-\delta}(\log t)^2\d t\ll_\epsilon T^{1/2+\epsilon},
\end{eqnarray*}
and
\begin{eqnarray*}
\int_{-\delta+iT}^{1+\delta+iT}\Delta(1-s){\zeta'\over \zeta}(s)\d s\ll \int_{-\delta}^{1+\delta} T^{-1/2+\sigma}(\log t)^2\d t\ll_\epsilon T^{1/2+\epsilon}
\end{eqnarray*}
by choosing $\delta=\epsilon/2$.

For the residue at a nontrivial $\zeta$-zero $\rho=\beta+i\gamma$ we find
\begin{eqnarray*}
{\rm{res}}_{s=\rho}=\lim_{s\to \rho}(s-\rho){\zeta'\over \zeta}(s)\Delta(1-s)=\Delta(1-\rho).
\end{eqnarray*}
Taking into account the functional equation, with $\rho$ also $1-\rho$ is a nontrivial $\zeta$-zero and $1-\overline{\rho}$ as well; since the sum over $\Lambda(n)$ is real and $\Delta$ on the real axis, there is no effect of conjugation. Hence, we arrive at  
\begin{eqnarray*}
\Sigma=\sum_{0<\gamma<T}\Delta(\rho).
\end{eqnarray*}
Substituting this together with the estimates above leads to the unconditional asymptotic formula 
\begin{eqnarray*}
-\psi(x)=\sum_{0<\gamma<T}\Delta(\rho)+O_\epsilon(T^{1/2+\epsilon}).
\end{eqnarray*}
In view of the equivalence of the Riemann hypothesis to (\ref{rh}) this proves Theorem \ref{2} for $a=0$.
\medskip

The initial idea for this paper actually was to find a new proof of the prime number theorem. For this purpose we applied Gonek's lemma (instead of Perron's formula) to the logarithmic derivative of $\zeta$. This is indeed the beginning of the proof above. 

Interestingly, the case $a\neq 0$ of Theorem \ref{2} needs a different reasoning (to our knowledge).

\section{Proof of Theorem \ref{2} in case $a\neq 0$}

Suppose that $a\neq 0$. Again we may assume that $T$ is at a distance from any $a$-point; more precisely, $\vert T-\gamma_a\vert\gg (\log T)^{-1}$, this time by (\ref{rvma})); this restriction can be removed at the expense of an error term of size $O(T^{1/2}\log T)$ by \eqref{rvma} and (\ref{delt}). 

By Cauchy's theorem, we have
\begin{eqnarray*}
\sum_{0<\gamma_a<T}\Delta(\rho_a)&=&{1\over 2\pi i}\int_R{\zeta'(s)\over \zeta(s)-a}\Delta(s)\d s,
\end{eqnarray*} 
where $R$ is the positive oriented rectangular contour with vertices $\kappa+i, \kappa+iT,-\delta+iT,-\delta+i$ with some small $\delta>0$ and a sufficiently large $\kappa>1$ such that there are no $a$-points to the right of $\kappa-\epsilon+i\R$, that is $\kappa=b_a+2\epsilon$, which is possible since Dirichlet series have a right half-plane free of zeros and $a$-points; the left vertical line segment does not cause any trouble since all but finitely many $a$-points in the left half-plane are distant from the line of integration (as already mentioned).

We split the integral into four; the two horizontals ones can be treated as in the case $a=0$ above with the estimate for $\zeta'(s)/(\zeta(s)-a)$. For the vertical intervals we proceed also quite similarly.

On the left we use the lower bound
\begin{equation*}
\zeta(\sigma+it)\gg_\epsilon t^{1/2-\sigma-\epsilon},
\end{equation*}
valid for $t\geq 2$ and $\sigma\leq0$, which follows easily from (\ref{delt}) and the functional equation (\ref{feq}) (or see \cite{gs}, Lemma 4). 
This allows to develop the logarithmic derivative into a convergent geometric series, namely
\begin{equation*}
{\zeta'(s)\over \zeta(s)-a}={\zeta'\over \zeta}(s)\cdot {1\over 1-a/\zeta(s)}={\zeta'\over \zeta}(s)\cdot \left\{1+{a\over \zeta(s)}+\sum_{k\geq 2}\left({a\over \zeta(s)}\right)^k\right\}.
\end{equation*}
This leads to
\begin{eqnarray*}
\int_{-\delta+iT}^{-\delta+i}{\zeta'(s)\over \zeta(s)-a}\Delta(s)\d s=I_1+I_2+I_3,
\end{eqnarray*} 
say, where
\begin{eqnarray*}
I_3=\int_{-\delta+iT}^{-\delta+i}{\zeta'\over \zeta}(s)\sum_{k\geq 2}\left({a\over \zeta(s)}\right)^k\Delta(s)\d s
\ll \int_1^T(\log t)^2 t^{2\epsilon-\delta-1/2}\d t\ll T^{1/2+2\epsilon}.
\end{eqnarray*} 
We continue with
\begin{eqnarray*}
I_1&=&\int_{-\delta+iT}^{-\delta+i}{\zeta'\over \zeta}(s)\Delta(s)\d s
\end{eqnarray*}
which we may rewrite by using the functional equation as
\begin{eqnarray*}
I_1&=&\int_{1+\delta-i}^{1+\delta-iT}\left({\Delta'\over \Delta}(1-s)-{\zeta'\over \zeta}(s)\right)\Delta(1-s)\d s=J_1+J_2,
\end{eqnarray*}
say. Using again (\ref{delt}), we find
\begin{eqnarray*}
J_1&=&-\int_{1+\delta-i}^{1+\delta-iT}\big( \Delta(1-s)\big)'\d s\ll T^{1/2+\delta}. 
\end{eqnarray*}
Moreover, $J_2$ turns out to be the conjugate of
\begin{eqnarray*}
\overline{J_2}&=&\int_{1+\delta+i}^{1+\delta+iT}-{\zeta'\over \zeta}(s)\Delta(1-s)\d s= 2\pi i\sum_{n\leq T/2\pi}\Lambda(n)+O_\delta(T^{1/2+2\delta}), 
\end{eqnarray*}
which, after applying the prime number theorem, contributes to the main term (as in the case $a=0$ above).

The integral $I_2$ depends on $a$; here we get
\begin{eqnarray*}
I_2&=&a\int_{1+\delta-i}^{1+\delta-iT}\left({\Delta'\over \Delta}(1-s)-{\zeta'\over \zeta}(s)\right){1\over \zeta(s)}\d s=H_1+H_2,
\end{eqnarray*}
say. The main term arises from $H_1$; to see that we apply the first derivative test (see, for example, \cite[Lemma 2.1]{ivic}) and get  
\begin{eqnarray*}
H_1&=&ia\int_1^T\big(\log {t\over 2\pi}+O(t^{-1})\big)\left(1+\sum_{n\geq 2}\mu(n)n^{-1-\delta-it}\right)\d t\\
&=&iaT\log{T\over 2\pi e}+O(\log T).
\end{eqnarray*}
In a similar way one bounds
\begin{eqnarray*}
H_2&=&a\int_{1+\delta-i}^{1+\delta-iT}\sum_{m\geq 2}\Lambda(m)m^{-s}\sum_{n\geq 1}\mu(n)n^{-s}\d s\ll 1.
\end{eqnarray*}
Thus,
\begin{eqnarray*}
I_2&=&iaT\log{T\over 2\pi e}+O(\log T).
\end{eqnarray*}

It remains to evaluate the vertical integral on the right, that is
\begin{eqnarray*}
\int_{\kappa+i}^{\kappa+iT}{\zeta'(s)\over \zeta(s)-a}\Delta(s)\d s.
\end{eqnarray*}
Using (\ref{partfrac}) and (\ref{delt}) this integral is bounded from above by $T^{3/2-\kappa}$. 

Collecting all results together, recalling that $\kappa>1$ and setting $\delta=\epsilon/2$, we obtain that
\begin{eqnarray*}
\sum_{0<\gamma_a<T}\Delta(\rho_a)&=&a{T\over 2\pi}\log{T\over 2\pi e}-\psi\left({T\over 2\pi}\right)+O_\epsilon\Big(T^{1/2+\epsilon}\Big).
\end{eqnarray*}
This implies Theorem \ref{2} in case $a\neq 0$.

\section{Proof of Theorem \ref{1} and a research question}

The statement follows from the distribution of values in the half-plane of absolute convergence of $\zeta$ and Lemma \ref{boj} of Bohr \& Jessen in particular. 

To see that, let $\epsilon>0$ be given. 
We assume that 
$$
\sum_{0<\gamma_a<T}\Delta(1-\rho_a)=o_\epsilon( T^{b_a-1/2-\epsilon}),
$$
as $T\to\infty$. By Lemma \ref{boj}, there exist infinitely many $a$-points $\rho_a^+=\beta_a^++i\gamma_a^+$ with $\vert \beta_a^+-b_a\vert<\epsilon$ and $\gamma_a^+\to\infty$.
It follows then that
$$
|\Delta(1-\rho_a^+)|=\left(\frac{\gamma_a^+}{2\pi}\right)^{\beta_a^+-1/2}\left(1+O\left(\frac{1}{\gamma_a^+}\right)\right)\geq\frac{1}{2} \left(\frac{\gamma_a^+}{2\pi}\right)^{b_a-1/2-\epsilon}
$$
for any sufficiently large $\gamma_a^+$.
Hence, if we denote by $\gamma_a'$ the term succeeding $\gamma_a^+$ in the increasing sequence $(\gamma_a)_{\gamma_a>0}$, we can obtain from the above relations that
\begin{align*}
\left|\sum_{0<\gamma_a<\gamma_a'}\Delta(1-\rho_a)\right|
\geq|\Delta(1-\rho_a^+)|-\left|\sum_{0<\gamma_a<\gamma_a^+}\Delta(1-\rho_a)\right|\gg_\epsilon\left({\gamma_a^+}\right)^{b_a-1/2-\epsilon}.
\end{align*}
It is not difficult to see that $\gamma_a'=\gamma_a^++O(1)$, $\gamma_a^+\to\infty$, or even $o(1)$, since the methods for obtaining such results are in their majority similar to the results regarding gaps of $\zeta$-zeros (see for example \cite[Chapter IX]{titch}). 
Therefore,
$$
\left({\gamma_a'}\right)^{b_a-1/2-\epsilon}\ll_\epsilon\left|\sum_{0<\gamma_a<\gamma_a^{'}}\Delta(1-\rho_a)\right|+o(1)=o_\epsilon\big(({\gamma_a'})^{b_a-1/2-\epsilon}\big),
$$
as $\gamma_a'\to\infty$, which is a contradiction. This finishes the proof of Theorem \ref{1}. 
\medskip

One could be tempted to improve upon this result by more advanced methods. We assume that $a\neq 1$ as well. Again we may rewrite the sum in question as a contour integral,
$$
\sum_{0<\gamma_a<T}\Delta(1-\rho_a)={1\over 2\pi i}\int_{R}{\zeta'(s)\over \zeta(s)-a}\Delta(1-s)\d s,
$$ 
where $R$ is a rectangular contour, similar to the one in the previous section with the difference that the right vertical segment has real part $\overline{\kappa}=\overline{\kappa}(\epsilon):=\overline{b_a}+\epsilon$.

The horizontal line segments can be treated as before; they can be bounded by $O_\epsilon\left(T^{\overline{\kappa}-1/2}\right)$ which is larger than their contributions in the previous cases.

For the integral over the left vertical line we have that
$$
{1\over 2\pi i}\int_{-\delta+iT}^{-\delta+i}{\zeta'(s)\over \zeta(s)-a}\Delta(1-s)\d s\ll\int_{1}^T(\log t)^2t^{-1/2+\delta}\d t\ll_\epsilon T^{1/2+\epsilon},
$$
where we set $\delta=\epsilon/2$.

For the integral over the right vertical line, we apply Lemma \ref{gonek} and obtain
\begin{equation}\label{ciurl}
{1\over 2\pi i}\int_{\overline{\kappa}+i}^{\overline{\kappa}+iT}{\zeta'(s)\over \zeta(s)-a}\Delta(1-s)\d s=\sum_{n\leq T/(2\pi)}\Lambda_a(n)+O_\epsilon\left(T^{\overline{\kappa}-1/2}\right). 
\end{equation}

Therefore,
$$
\sum_{0<\gamma_a<T}\Delta(1-\rho_a)=\sum_{n\leq T/(2\pi)}\Lambda_a(n)+O_\epsilon\left(T^{\overline{\kappa}-1/2}\right).
$$
This is obviously a more general formula then the one in Theorem 1 for the case of $a=0$ since $\Lambda_0=-\Lambda$ and $\overline{b_0}=1$.
Also, for any $a\neq1$ such that $\overline{b_a}-b_a<1/2$, the formula above and Theorem 2 imply that
$$
\Omega_\epsilon\left(T^{b_a-1/2-\epsilon}\right)=\sum_{0<\gamma_a<T}\Delta(1-\rho_a)\ll_\epsilon T^{b_a+\epsilon}.
$$

Applying Perron's formula (see \cite[\S 3.12]{titch}), we have
$$
\sum_{n\leq x}\Lambda_a(n)={1\over 2\pi i}\int_{\overline{\kappa}-i\tau}^{\overline{\kappa}+i\tau}{\zeta'(s)\over \zeta(s)-a}{x^s\over s}\d s+O_\epsilon\left({x^{\kappa+1}\over \tau}\right),
$$
where $\tau>0$ is a parameter and $x=T/(2\pi)$ or any other positive quantity having a fixed distance from integers. 
Shifting the line of integration to the left, we get contributions from the residues of the integrand at the $a$-points $\rho_a=\beta_a+i\gamma_a$ of $\zeta$ as well as from the simple pole of $\zeta(s)$ at $s=1$. The latter residue, however, turns out to be small compared to
$$
{\rm{res}}_{s=\rho_a}=\lim_{s\to\rho_a}(s-\rho_a){\zeta'(s)\over \zeta(s)-a}{x^s\over x}={x^{\rho_a}\over \rho_a},
$$
which is of size $T^{\beta_a}/\gamma_a$ in our situation.

 If there is not too much cancelation for the sum of these residues, one could expect their contribution to be of order $T^{b_a}\sum_{0<\gamma<T}\gamma_a^{-1}$ which, in view of Lemma \ref{boj}, is of size $T^{b_a}\log T$. Actually, a hypothetical $\zeta$-zero off the critical line (and also if there are more than just one), would have a similar effect on the error term in the prime number theorem (see \cite[Chapter V]{ingham}).

We have imposed above the condition $\overline{b_a}-b_a<1/2$ to deduce an upper bound for the sum of $\Delta(1-\rho_a)$. It is natural to ask for which complex numbers $a$ this condition holds, when in general we know that the maximum difference of the abscissa of convergence and absolute convergence can be $1$ (see \cite[\S 9.13]{tit}).

A straightforward way to describe such points is to bound $\overline{b_a}$ by $3/2$ since $b_a>1$ for $a\neq1$. Let $\sigma>1$ be sufficiently large such that $|\zeta(s)-1|=2^{-\sigma}|(1+o(1))|<|a-1|$. Then
\begin{eqnarray*}
\frac{1}{\zeta(s)-a}&=&\frac{1}{1-a}\left(1-\frac{\zeta(s)-1}{a-1}\right)^{-1}\\
&=&\frac{1}{1-a}\sum_{k\geq 0}\left(\frac{\zeta(s)-1}{a-1}\right)^{k}=-\sum_{k\geq 0}\sum_{n\geq1}\frac{d_k^*(n)}{{(a-1)^{k+1}}n^s},
\end{eqnarray*}
where $d^*_k(n)$ denotes the number of decompositions of $n$ into a product of $k$ factors greater than $1$. We observe that
$$
\sum_{k\geq 0}\sum_{n\geq1}\left|\frac{d_k^*(n)}{{(a-1)^{k+1}}n^s}\right|=\frac{1}{|a-1|}\sum_{k\geq 0}\sum_{n\geq1}\frac{d_k^*(n)}{{|a-1|^{k}}n^\sigma}=\frac{1}{|a-1|-|\zeta(\sigma)-1|}
$$
for any $\sigma>\sigma^*$ where $\zeta(\sigma^*)-1=|\zeta(\sigma^*)-1|=|a-1|$.

Thus, for any such $\sigma$ the double series converges absolutely and we can interchange summation to obtain a Dirichlet series representation of the reciprocal of $\zeta(s)-a$ and, consequently, of its logarithmic derivative for which we are going to have that $\overline{b_a}\leq\sigma^*$. Hence, {\it for every complex number $a$ satisfying $\vert a-1\vert>\zeta(3/2)-1\approx1.6123...$ it follows that $\overline{b_a}<3/2$.} On the other hand, for such $a$ the Omega-theorem does not imply that the non-existence of the mean of $\Delta(1-\rho_a)$. One can also see from the construction above that if $a>1$ then $b_a=\overline{b_a}=\sigma^*$ which also serves our purpose.

\section*{Acknowledgement}
The third author was supported by JSPS KAKENHI Grant Numbers 18K13400 and 22K13895, and also MEXT Initiative for Realizing Diversity in the Research Environment.

\bigskip

\noindent
Athanasios Sourmelidis\\
Institute of Analysis and Number Theory, TU Graz\\
Steyrergasse 30, 8010, Graz, Austria\\
sourmelidis@math.tugraz.at
\bigskip

\noindent
J\"orn Steuding\\
Department of Mathematics, W\"urzburg University\\
Emil Fischer-Str. 40, 97\,074 W\"urzburg, Germany\\
steuding@mathematik.uni-wuerzburg.de
\bigskip

\noindent
Ade Irma Suriajaya\\
Faculty of Mathematics, Kyushu University \\
744 Motooka, Nishi-ku, Fukuoka 819-0395, Japan\\
adeirmasuriajaya@math.kyushu-u.ac.jp

\end{document}